\documentclass[12pt]{article}
\usepackage{amsfonts}
\usepackage{amssymb}
\usepackage{colortbl}
\usepackage{color}

\newtheorem{theorem}{Theorem}[section]

\newtheorem{definition}[theorem]{Definition}

\newtheorem{lemma}[theorem]{Lemma}

\begin{document}

\renewcommand{\theequation}{\arabic{section}.\arabic{equation}}

\bigskip

\begin{center}
{\Large Approximating technique for the solution to the }${\LARGE H}^{\infty
}${\Large -optimal control problem}

\bigskip

Gabriela Marinoschi

\medskip

\textquotedblleft Gheorghe Mihoc-Caius Iacob\textquotedblright\ Institute of
Mathematical Statistics and

Applied Mathematics of the Romanian Academy,

Calea 13 Septembrie 13, Bucharest, Romania

E-mail: gabriela.marinoschi@acad.ro

\bigskip
\end{center}

\noindent \noindent Abstract. This paper continues the results published in 
\cite{GM-2026}, by presenting an approximating method of calculating the
solution to the $H^{\infty }$-optimal control problem, which solves the
robust stabilization of a system by a feedback controller, and to the
associated Riccati equation. This method generates a sequence of infinite
dimensional systems approximating the original system and the solution to
the approximating $H^{\infty }$-optimal control problem is proved to
converge to the solution of the original $H^{\infty }$-optimal control
problem. An example for calculating the solution to the Riccati equation and
the feedback controller is provided at the end.

\medskip

Keywords: $H^{\infty }$-control, feedback control, robust stability, Riccati
equation, approximation method

\medskip

MSC 2020: 93B36, 93B52, 93D09

\section{Introduction}

\setcounter{equation}{0}

The $H^{\infty }$-optimal control is a technique used in control theory to
design robust stabilizing feedback controllers that force a system to
achieve stability with a prescribed performance despite the uncertainties,
such as parameter variations, external disturbances, observation sensor
noise, which may corrupt the reliable functionality of the system. The
feedback control determined via the $H^{\infty }$-optimal control technique
for stabilizing the perturbed system has an expression constructed by the
aid of the solution to an associated algebraic Riccati equation. The origins
of the $H^{\infty }$-optimal control technique developed for
finite-dimensional systems are found in the papers \cite{Glover}, \cite%
{Zames}, continued then by other papers e.g., \cite{Ichikawa-91}, \cite%
{Zhou-Doyle-Glover}. The extension to the infinite-dimensional systems with
bounded control operators was done in \cite{VB-hyp-92} and \cite{VB-hyp-95}
and focused on the $H^{\infty }$-boundary control under state feedback for
hyperbolic-like systems, and also in \cite{Bensoussan}, \cite{Ichikawa-96},
and in the monograph \cite{Keulen} addressing the $H^{\infty }$-control
problem for a class of systems with unbounded control operators. Let us
recall that in the finite-dimensional case the feedback operator is found by
solving a matrix algebraic Riccati equation. In the case of
infinite-dimensional systems this can be done by solving an
infinite-dimensional operator algebraic Riccati equation. For a general
theory regarding the Riccati equation in the infinite-dimensional case and
its applications we refer the reader to the monographs \cite{IL-RT-1991} and 
\cite{IL-RT-1999}. An operatorial method for solving the Riccati equation $%
A^{\ast }P+PA+P^{2}=F$ was given in \cite{VB-76}, Chapter II, 3.3, see also 
\cite{Brezis-1970}. A study for a very specific heat equation was done in 
\cite{BMR-2020} and an approximation theory of solutions to operator Riccati
equations for $H^{\infty }$- optimal control have been developed e.g., in 
\cite{Ito} and \cite{Morris} for bounded control operators. The $H^{\infty }$%
-optimal control abstract problem has been solved in \cite{GM-23-Hinf} for
linear infinite dimensional systems of parabolic type, relying on a
variational technique. In what concerns the determination of the feedback
operator by numerical techniques, we can cite e.g., the papers \cite%
{Har-Kalm}, where the problem was formulated as a nonlinear operator
equation set up in a sparse grid discretization, and \cite{Guo-Tan} treating
the Dirichlet boundary control by employing the Galerkin approximation
generating a sequence of solvable finite-dimensional systems that
approximate the original infinite-dimensional system and demonstrating that
the sequence of solutions to the associated finite-dimensional algebraic
Riccati equations converges in norm to the solution of the
infinite-dimensional operator algebraic Riccati equation. The paper \cite%
{GM-2026} deals with the existence of the solution to the associated
algebraic Riccati equation by a direct operatorial approach in the space of
Hilbert-Schmidt operators. The proofs develop a constructive approach, by
calculating the solution by a Fourier series along a suitable orthonormal
base of the Hilbert space of the leading operator in the system to be
stabilized. However, in a concrete situation the calculation of an infinite
series is not possible, so that the current paper continues the results of
the paper \cite{GM-2026}, by presenting an approximating method of replacing
the feedback operator by a truncated series up to a finite number of terms,
in order to allow its numeric calculation.

The outline of the paper is the following. At the beginning there are
presented some preliminaries introducing the $H^{\infty }$-optimal control
problem and some previous results concerning its solution, as well as the
solution to the associated algebraic Riccati equation. Then, the Fourier
series expansion representing the optimal control found for the $H^{\infty }$%
-optimal control problem is truncated up to a certain number of terms and
the the corresponding approximating equation is studied for obtaining
results of existence, uniqueness and some estimates in Theorem \ref{T2.2}.
Moreover, this approximating control is proved to be a solution, in some
sense, to the approximated $H^{\infty }$-optimal control problem. Relying on
estimates it is proved still in Theorem \ref{T2.2}, the the limit of the
approximating solution tends to the solution of the original $H^{\infty }$%
-optimal control problem. All these are proved first in the case when the
main operator in the system is coercive. The singular case when the operator
is only positive is discussed in Section 3. Another approximation method is
described in the case when the operators involved in the system to be
stabilized are not Hilbert-Schmidt. The approach considered here has the
advantage of putting in evidence constructive proofs which set up a method
for calculating the numerical solution in concrete cases. An example showing
the way in which the theoretical results can be applied is provided in
Section 5.

\subsection{Notation and some definitions}

We put down a few notation, definitions and results necessary for
approaching the problem. Let $X$ be a real Hilbert space with the dual $%
X^{\prime }.$ By $(\cdot ,\cdot )_{X}$, $\left\Vert \cdot \right\Vert _{X}$
and $\left\langle \cdot ,\cdot \right\rangle _{X^{\prime },X}$ we denote the
scalar product, the norm and the pairing between $X^{\prime }$ and $X,$
respectively. If $Y$ is another Hilbert space, $L(X,Y)$ represents the space
of all linear continuous operators from $X$ to $Y.$ If $F\in L(X,X^{\prime
}) $ the adjoint of $F$ is denoted by $F^{\ast }.$ Moreover, let us consider
the Hilbert spaces $V,$ $H,$ $V^{\prime }$ forming a variational triplet, $%
V\subset H\subset V^{\prime }$ with dense and continuous injections and $V$
compact in $H.$

\noindent We also recall that a linear operator, $L:V\rightarrow V^{\prime }$
has the property $L\geq 0$ if $\left\langle Ly,y\right\rangle _{V^{\prime
},V}\geq 0$ for all $y\in V,$ and $L>0$ if $\left\langle Ly,y\right\rangle
_{V^{\prime },V}>0,$ for all $y\in V,$ $y\neq 0.$ The pairing reduces to the
scalar product in $H$ if $L\in L(H,H).$

\noindent Now, we introduce the Hilbert-Schmidt spaces that represent the
most appropriate functional spaces for solving the operatorial algebraic
Riccati equation, as well as the associate $H^{\infty }$-optimal control
problem.

\noindent Let $H_{i},$ $i=1,2$ be separable Hilbert spaces. The operator $%
P:H_{1}\rightarrow H_{2}$ is said to be Hilbert-Schmidt if for some complete
orthonormal basis $\{e_{n}\}_{n\geq 1}$ of $H_{1}$ we have $%
\sum\limits_{n=1}^{\infty }(Pe_{n},Pe_{n})_{H_{2}}<\infty .$

\noindent We denote by $(H_{1},H_{2})_{H.S.}$ the space of Hilbert-Schmidt
operators from $H_{1}$ to $H_{2}$ and recall that it is a Hilbert space with
the scalar product and norm%
\[
(P,Q)_{(H_{1},H_{2})_{H.S.}}=\sum\limits_{n=1}^{\infty
}(Pe_{n},Qe_{n})_{H_{2}},\mbox{ } 
\]%
\begin{equation}
\left\Vert P\right\Vert
_{(H_{1},H_{2})_{H.S.}}^{2}=\sum\limits_{n=1}^{\infty }\left\Vert
Pe_{n}\right\Vert _{H_{2}}^{2},\mbox{ }P,\mbox{ }Q\in (H_{1},H_{2})_{H.S.}.
\label{11-0}
\end{equation}%
We denote 
\begin{equation}
\mathcal{H}:=(H,H)_{H.S.},\mbox{ }\mathcal{V}:=(V^{\prime },H)_{H.S.}\cap
(H,V)_{H.S.}.  \label{11-2}
\end{equation}%
The dual space of $\mathcal{V},$ denoted $\mathcal{V}^{\prime }$ can be
identified with $(V,H)_{H.S.}+(H,V^{\prime })_{H.S.}$ and we have the
embedding $\mathcal{V}\subset \mathcal{H}\subset \mathcal{V}^{\prime }$ with
dense and continuous inclusions.

\subsection{$H^{\infty }$-control problem and the Riccati equation}

We begin by briefly explaining the meaning of the $H^{\infty }$-control
problem for the linear system%
\begin{equation}
y^{\prime }(t)+Ay(t)=B_{1}w(t)+B_{2}u(t),\mbox{ }t>0,  \label{1}
\end{equation}%
\begin{equation}
z(t)=C_{1}y(t)+D_{1}u(t)\mbox{ }t>0,  \label{2}
\end{equation}%
\[
y(0)=y_{0}, 
\]%
where $y$ is the state, $w$ is a unknown perturbation, $u$ is the control to
be found and $z$ is the output. For $C_{1}$ and $D_{1}$ we assume, as
technical simplifications, the conditions $D_{1}^{\ast }C_{1}=0$ and $%
D_{1}D_{1}^{\ast }=1,$ considered in \cite{GM-23-Hinf}, as well as $D_{1}=0$
leading to the generally used output form $z(t)=C_{1}y(t).$

Let $H,$ $U,$ $W,$ $Z$ be real Hilbert spaces, identified with their duals.
The hypotheses assumed in the paper are:

$(i_{1})$ $A\in L(V,V^{\prime })$, $\left\langle Ay,y\right\rangle
_{V^{\prime },V}\geq \omega \left\Vert y\right\Vert _{V}^{2},$ for all $y\in
V,$ $\omega \geq 0,$

$(i_{2})$ $B_{1}\in L(W,H),$ $B_{2}\in L(U,H),$ $C_{1}\in L(H,Z),$ $D_{1}\in
L(U,Z),$

\noindent where $W$ is a space for the perturbation, $Z$ a space for the
output and $U$ the space of control.

The $H^{\infty }$-optimal control technique means to design a feedback
control $u=Fy$ that stabilizes exponentially the system and forces it to
achieve a robust performance, that is to maintain performance and stability
under external disturbance, ensuring thus that the system operates reliably
under different conditions. More exact, the problem reduces to determine $%
F:H\rightarrow U$ such that the operator $-A+B_{2}F$ generates an
exponentially stable semigroup on $H,$ that is 
\begin{equation}
\left\Vert e^{(-A+B_{2}F)t}y\right\Vert _{H}\leq Ce^{-\alpha t}\left\Vert
y\right\Vert _{H}\mbox{ for all }y\in H,\mbox{ with }C>0,\mbox{ }\alpha >0,
\label{2-0}
\end{equation}%
and the output corresponding to the controlled system, that is with $u=Fy,$
has the property 
\begin{equation}
\left\Vert z\right\Vert _{L^{2}(\mathbb{R}_{+};Z)}=\left\Vert
C_{1}y+D_{1}u\right\Vert _{L^{2}(\mathbb{R}_{+};Z)}\leq \gamma \left\Vert
w\right\Vert _{L^{2}(\mathbb{R}_{+};W)}  \label{3}
\end{equation}%
for $\gamma $ a prescribed positive value called the robust performance
index. In the latter relation $z$ is viewed as depending on $w$ by solving (%
\ref{1}) with $u=Fy.$ We recall that by the Datko's result (see \cite%
{Datko-72}), relation (\ref{2-0}) is equivalent with $e^{(-A+B_{2}F)t}y\in
L^{2}(0,\infty ;H).$

The robust stabilization via the $H^{\infty }$-optimal control problem
expressed by the relations (\ref{2-0}) and (\ref{3}) is solved by the
results of Theorems 2.1 and 4.1 in \cite{GM-2026}, concentrated here in a
single statement:

\begin{theorem}
\label{T1.1}Let $\gamma >0$ and let us denote 
\begin{equation}
\Gamma :=B_{2}B_{2}^{\ast }-\gamma ^{-2}B_{1}B_{1}^{\ast }.  \label{12-0}
\end{equation}%
In addition to $(i_{1})$ and $(i_{2})$ with $\omega >0,$ let 
\begin{equation}
C_{1}\in \mathcal{H},\mbox{ }\Gamma \in \mathcal{H},\mbox{ }\Gamma \geq 0.
\label{12}
\end{equation}%
Then, the Riccati equation 
\begin{equation}
A^{\ast }P+PA+P(B_{2}B_{2}^{\ast }-\gamma ^{-2}B_{1}B_{1}^{\ast
})P=C_{1}^{\ast }C_{1}  \label{5}
\end{equation}%
has a unique solution $P\in L(H,H),$ with 
\begin{equation}
P\geq 0,\mbox{ }P=P^{\ast }.  \label{15}
\end{equation}%
Assume still that $w\in L^{2}(0,\infty ;W)$ and $\Lambda
=-A-B_{2}B_{2}^{\ast }P$ generates an exponentially stable semigroup on $H.$
Then, the feedback operator 
\begin{equation}
F=-B_{2}^{\ast }P  \label{4}
\end{equation}%
solves the $H^{\infty }$-control problem associated to system (\ref{1})-(\ref%
{2}), that is relation (\ref{3}).
\end{theorem}

The Riccati equation (\ref{5}) has a solution also in the singular case when 
$A\geq 0$ (see \cite{GM-2026}, Theorem 3.1) and so we have

\begin{theorem}
\label{T1.2}Assume $(i_{1}),$ $(i_{2})$ with $\omega =0,$ and 
\begin{equation}
C_{1}\in (V^{\prime },H)_{H.S.},\mbox{ }\ker C_{1}=\{0\},  \label{51}
\end{equation}%
\begin{equation}
\Gamma \in \mathcal{H},\mbox{ }(\Gamma y,y)_{H}\geq g_{0}\left\Vert
y\right\Vert _{H}^{2},\mbox{ for any }y\in H,\mbox{ }g_{0}>0.  \label{52}
\end{equation}%
Then, there exists a unique $P\in \mathcal{V}$ solution to the algebraic
Riccati equation (\ref{5}), having the properties (\ref{15}). The second
part of the Theorem \ref{T1.1} holds true.
\end{theorem}

The main step in approaching the $H^{\infty }$-optimal control problem is
thus to solve the algebraic Riccati equation (\ref{5}). The previous results
will be used as starting points in achieving the aim of this paper.

\section{The coercive case}

\setcounter{equation}{0}

We consider here the coercive case for the operator $A,$ that is $(i_{1})$
with $\omega >0.$ By Theorem \ref{T1.1} the solution $P$ to the Riccati
equation (\ref{5}) exists and it is unique and there exists a feedback
operator $F=-B_{2}^{\ast }P$ which solves the $H^{\infty }$-optimal control
problem, correspondingly relation (\ref{3}), providing the control $%
u=-B_{2}^{\ast }Py.$

Moreover, we assert that the solution $P$ can be constructed by an expansion
in Fourier series along an orthonormal basis $\{e_{j}\}_{j\geq 1}$ in $H.$
Further, we briefly explain this. For convenience, the basis can be
considered formed by the eigenvectors generated by the canonical isomorphism 
$J:V\rightarrow V^{\prime },$ $Je_{j}=\rho _{j}^{2}e_{j},$ where $\{\rho
_{j}^{2}\}_{j\geq 1}$ is the sequence of the positive eigenvalues
increasingly ordered. Then, $y\in H$ is expressed in Fourier series as 
\begin{equation}
y=\sum_{j=1}^{\infty }y_{j}e_{j},\mbox{ for all }y\in H,\mbox{ }%
y_{j}=(y,e_{j})_{H},\mbox{ }j\geq 1.  \label{52-0}
\end{equation}%
For $P,$ the solution to (\ref{5}), we choose the expansion 
\begin{equation}
Py=\sum_{j=1}^{\infty }\gamma _{j}y_{j}e_{j},\mbox{ for
all }\mbox{ }y\in H,  \label{52-1}
\end{equation}%
implying that 
\begin{equation}
Pe_{j}=\gamma _{j}e_{j},\mbox{ }j\geq 1,  \label{52-3}
\end{equation}%
where $\gamma _{j}$ should be positive. Indeed, it can be proved that there
exist $\gamma _{j}\geq 0$ such that $P$ defined in this way verifies
equation (\ref{5}). To this end, we replace the previous expressions in (\ref%
{5}) applied to $y=e_{j}$ and get after a few calculations, since $\Gamma
=\Gamma ^{\ast },$ that 
\[
(A^{\ast }P+PA)e_{j}+\sum_{k=1}^{\infty }\gamma _{j}\gamma _{k}e_{k}(\Gamma
e_{j},e_{k})_{H}=C_{1}^{\ast }C_{1}e_{j},\mbox{ }j\geq 1. 
\]%
This multiplied scalarly by $e_{j}$ yields 
\[
2\left\langle Ae_{j},Pe_{j}\right\rangle _{V^{\prime },V}+(\gamma
_{j})^{2}(\Gamma e_{j},e_{j})_{H}=(C_{1}^{\ast }C_{1}e_{j},e_{j})_{H},%
\mbox{
}j\geq 1. 
\]%
Using again (\ref{52-3}) we obtain the equation 
\[
(\gamma _{j})^{2}g_{jj}+2\gamma _{j}a_{jj}-\left\Vert C_{1}e_{j}\right\Vert
_{H}^{2}=0,\mbox{ }j\geq 1, 
\]%
where $g_{jj}:=(\Gamma e_{j},e_{j})_{H}$ and $a_{jj}:=\left\langle
Ae_{j},e_{j}\right\rangle _{V^{\prime },V}.$ We recall that by (\ref{12}), $%
(\Gamma y,y)\geq 0$ and so $g_{jj}\geq 0.$ Also, by $(i_{1})$ with $\omega
>0 $ it follows that $a_{jj}\geq \omega >0.$ The above equation has a
nonnegative solution%
\begin{equation}
\gamma _{j}=\left\{ 
\begin{tabular}{ll}
$\frac{-a_{jj}+\sqrt{(a_{jj})^{2}+g_{jj}\left\Vert C_{1}e_{j}\right\Vert
_{H}^{2}}}{g_{jj}}$ & $\mbox{if }g_{jj}\neq 0,\mbox{ }$ \\ 
$\frac{\left\Vert C_{1}e_{j}\right\Vert _{H}^{2}}{2a_{ij}},$ & if $g_{jj}=0$%
\end{tabular}%
\right.  \label{52-2}
\end{equation}%
so that the eigenvalues of $P$ in the basis $\{e_{j}\}_{j\geq 1}$ are found
to be positive. They are also bounded, $0\leq \gamma _{j}\leq \frac{%
\left\Vert C_{1}e_{j}\right\Vert _{H}^{2}}{2\omega },$ $j\geq 1.$ We can
check that $P$ constructed in this way belongs to $\mathcal{H},$ that is the
sum of the squares of $\gamma _{j}$ is finite. Indeed, since the series with
the terms $\left\Vert C_{1}e_{j}\right\Vert _{H}^{2}$ converges because its
limit is $\left\Vert C_{1}\right\Vert _{\mathcal{H}}^{2}$, we deduce that
the series $\sum_{j=1}^{\infty }\left\Vert C_{1}e_{j}\right\Vert _{H}^{4}$
converges too and this verifies that $\left\Vert P\right\Vert _{\mathcal{H}%
}^{2}=\sum_{j=1}^{\infty }\gamma _{j}^{2}<\infty .$

In conclusion, $Py$ can be calculated by (\ref{52-1}) with $\gamma _{j}$
given by (\ref{52-2}) and the control $u$ is given by $u=-B_{2}^{\ast }Py$
(by Theorem \ref{T1.1}). By plugging these in equation (\ref{1}), it becomes%
\begin{equation}
y^{\prime }(t)=-Ay(t)-B_{2}B_{2}^{\ast }Py(t)+B_{1}w(t)  \label{90}
\end{equation}%
and it is clear that the solution $y$ and consequently $u$ depend on $w.$

However, if numerical computations are necessary, the infinite series (\ref%
{52-1}) cannot be computed, so that one has to consider only an appropriate
finite number $N$ of terms in the sum. Thus, we consider a truncated
expression of the feedback control 
\begin{equation}
u_{N}=-B_{2}^{\ast }P_{N}y,\mbox{ }F_{N}=-B_{2}^{\ast }P_{N},  \label{13-0}
\end{equation}%
where%
\begin{equation}
P_{N}y=\sum_{j=1}^{N}\gamma _{j}y_{j}e_{j},\mbox{ for all }y\in H,
\label{13}
\end{equation}%
which is obviously positive and selfadjoint, and 
\begin{equation}
P_{N}e_{k}=\left\{ 
\begin{tabular}{ll}
$\gamma _{k}e_{k},$ & $k=1,...,N$ \\ 
$0,$ & $k\geq N+1.$%
\end{tabular}%
\right.  \label{13-00}
\end{equation}%
Now, let us introduce $F_{N}$ instead of $F=-B_{2}^{\ast }P$ in system (\ref%
{90}), getting 
\[
y^{\prime }(t)+Ay(t)=-B_{2}B_{2}^{\ast }P_{N}y(t)+B_{1}w(t),\mbox{ }t>0. 
\]

The first step is to show that this equation has a solution and since it
depends on $N,$ this entitles us to denote it $y_{N}.$ So, we have to study
the system%
\begin{equation}
y_{N}^{\prime }(t)+Ay_{N}(t)=-B_{2}B_{2}^{\ast }P_{N}y_{N}(t)+B_{1}w(t),%
\mbox{ }t>0,  \label{20}
\end{equation}%
\begin{equation}
z_{N}(t)=C_{1}y_{N}(t)+D_{1}u_{N}(t)\mbox{ }t>0,  \label{21}
\end{equation}%
with $y_{N}(0)=y_{0}\in H.$

Simultaneously, we shall show that the truncated control 
\begin{equation}
u_{N}=-B_{2}^{\ast }P_{N}y_{N},\mbox{ }  \label{11}
\end{equation}%
solves in some approximating sense the $H^{\infty }$-optimal control problem.

The second step is that, by relying on appropriate estimates for the
solution $y_{N},$ to pass to the limit in (\ref{20})-(\ref{11}) as $%
N\rightarrow \infty $ and to prove that $y_{N}$ and $u_{N}$ are suitable
approximations for the solution $y$ to (\ref{90}) and for the control $u,$
respectively.

We give now the definition of the approximating solution to the $H^{\infty }$%
-control problem related to system (\ref{20})-(\ref{11}).

\begin{definition}
\label{D2.1}Let $\gamma >0$ be fixed. The feedback control $%
u_{N}=-B_{2}^{\ast }P_{N}y_{N}$ is called an approximating solution to the $%
H^{\infty }$-control problem associated to system (\ref{20})-(\ref{11}) if
the operator $-A-B_{2}B_{2}^{\ast }P_{N}$ generates an exponentially stable
semigroup and if the output satisfies 
\begin{equation}
\left\Vert z_{N}\right\Vert _{L^{2}(\mathbb{R}_{+};Z)}=\left\Vert
C_{1}y_{N}+D_{1}u_{N}\right\Vert _{L^{2}(\mathbb{R}_{+};Z)}\leq \gamma
\left\Vert w\right\Vert _{L^{2}(\mathbb{R}_{+};W)}+C\delta _{N}  \label{22}
\end{equation}%
where $y_{N}$ is the solution to (\ref{20}) and $\{\delta _{N}\}_{N\geq 1}$
is a sequence depending on $N,$ such that $\delta _{N}\rightarrow 0$ as $%
N\rightarrow \infty .$
\end{definition}

\noindent By $C$ we denote further several positive constants which may
change from line to line.

\noindent The first result is the following:

\begin{theorem}
\label{T2.2}Let assumptions of Theorem 1.1 hold. Let $P$ be the solution to
the Riccati equation (\ref{5}) and $\varepsilon
_{N}:=\sum\limits_{j=N+1}^{\infty }\gamma _{j}^{2}$ with $\gamma _{j}$ the
eigenvalues of $P,$ by (\ref{52-3}). Then, equation (\ref{20}) has a unique
solution 
\begin{equation}
y_{N}\in C([0,\infty );H)\cap L^{2}(0,\infty ;H)  \label{23}
\end{equation}%
and the corresponding control $u_{N}=-B_{2}^{\ast }P_{N}y_{N}$ is the
approximating solution to the $H^{\infty }$-control problem associated to
system (\ref{20})-(\ref{11}). Moreover, the following estimates hold true:%
\begin{equation}
\left\Vert y_{N}\right\Vert _{C([0,T];H)\cap L^{2}(0,T;V)}\leq C,\mbox{ for
all }T>0,  \label{24}
\end{equation}%
\begin{equation}
\left\Vert (P-P_{N})y\right\Vert _{H}^{2}\leq \varepsilon _{N}\left\Vert
y\right\Vert _{H}^{2},\mbox{ for all }y\in H,  \label{25-1}
\end{equation}%
\begin{equation}
\left\Vert P_{N}\right\Vert _{\mathcal{H}}\leq C.  \label{25-0}
\end{equation}%
For $N$ large enough we have%
\begin{equation}
\int_{0}^{\infty }\left\Vert y_{N}(t)\right\Vert _{H}^{2}dt\leq C,
\label{25}
\end{equation}%
\begin{equation}
\int_{0}^{\infty }\left\Vert y_{N}(t)-y(t)\right\Vert _{H}^{2}dt\leq
C\varepsilon _{N},  \label{25-00}
\end{equation}%
\[
\int_{0}^{\infty }\left\Vert u_{N}(t)-u(t)\right\Vert _{H}^{2}dt\leq
C\varepsilon _{N}, 
\]%
\[
\int_{0}^{\infty }\left\Vert z_{N}(t)-z(t)\right\Vert _{H}^{2}dt\leq
C\varepsilon _{N}. 
\]
Finally, we have%
\[
y_{N}\rightarrow y\mbox{ strongly in }L^{2}(0,\infty ;H)\mbox{, as }%
N\rightarrow \infty , 
\]%
\[
y_{N}(t)\rightarrow y(t)\mbox{ strongly in }H,\mbox{ as }N\rightarrow \infty
,\mbox{ uniformly in }t>0, 
\]%
\[
P_{N}\rightarrow P\mbox{ strongly in }L(H,H)\mbox{ and weakly in }\mathcal{H}%
\mbox{ as }N\rightarrow \infty , 
\]%
\[
u_{N}\rightarrow u\mbox{ strongly in }L^{2}(0,\infty ;U)\mbox{ as }%
N\rightarrow \infty , 
\]%
\[
u_{N}(t)\rightarrow u(t)\mbox{ strongly in }U,\mbox{ a.e. }t>0,\mbox{ as }%
N\rightarrow \infty , 
\]%
\[
z_{N}\rightarrow z\mbox{ strongly in }L^{2}(0,\infty ;Z)\mbox{ as }%
N\rightarrow \infty , 
\]%
\[
z_{N}(t)\rightarrow z(t)\mbox{ strongly in }Z,\mbox{ a.e. }t>0,\mbox{ as }%
N\rightarrow \infty , 
\]%
where $y$ is the solution to (\ref{90}) and $u$ is the solution to the
original $H^{\infty }$-control problem.
\end{theorem}

\noindent \textbf{Proof. }Under assumptions of Theorem \ref{T1.1}, the%
\textbf{\ }$H^{\infty }$-control problem associated to system (\ref{1})-(\ref%
{2}) has the solution $u=-B_{2}^{\ast }Py$ with $y$ the solution to (\ref{90}%
) and $P$ the unique solution to the Riccati equation (\ref{5}). We specify
that $y\in C([0,\infty );H)\cap L^{2}(0,\infty ;H),$ since the operator $%
\Lambda $ generates an exponentially stable semigroup.

Passing to equation (\ref{20}), since $A$ is linear continuous and coercive
it is $m$-accretive (maximal monotone on $H$) and consequently $-A$
generates a $C_{0}$-semigroup of contractions in $H$. Then, since $%
B_{2}B_{2}^{\ast }P_{N}$ is linear continuous on $H$ it follows that $%
\Lambda _{N}:=-A-B_{2}B_{2}^{\ast }P_{N}$ generates a $C_{0}$-semigroup and
so equation (\ref{20}) has a unique solution $y_{N}\in C([0,\infty );H).$ We
are going first to prove that $\Lambda _{N}$ is exponentially stable.

At the beginning we establish two estimates. We see that $P_{N}\in \mathcal{H%
}$ since 
\begin{equation}
\left\Vert P_{N}\right\Vert _{\mathcal{H}}^{2}=\sum_{j=1}^{\infty
}\left\Vert P_{N}e_{j}\right\Vert _{H}^{2}=\sum_{j=1}^{N}\left\Vert \gamma
_{j}e_{j}\right\Vert _{H}^{2}\leq \left\Vert P\right\Vert _{\mathcal{H}}^{2}.
\label{30-0}
\end{equation}%
Then, 
\begin{equation}
\left\Vert P_{N}\widetilde{y}\right\Vert _{H}^{2}\leq \left\Vert
\sum_{j=1}^{N}\gamma _{j}e_{j}(\widetilde{y},e_{j})_{H}\right\Vert
_{H}^{2}\leq \left\Vert \widetilde{y}\right\Vert _{H}^{2}\left\Vert
\sum_{j=1}^{N}\gamma _{j}^{2}\right\Vert ^{2}\leq \left\Vert \widetilde{y}%
\right\Vert _{H}^{2}\left\Vert P\right\Vert _{\mathcal{H}}^{2},  \label{30}
\end{equation}%
for $\widetilde{y}\in H,$ where we took into account (\ref{52-3}) and (\ref%
{13-00}).

Next, let us write (\ref{20}) as 
\begin{equation}
y_{N}^{\prime }(t)=(-A-B_{2}B_{2}^{\ast }P)y_{N}(t)+(B_{2}B_{2}^{\ast
}P-B_{2}B_{2}^{\ast }P_{N})y_{N}(t)+B_{1}w(t),  \label{26}
\end{equation}%
and calculate some estimates, taking advantage of the fact that the operator 
$\Lambda =-A-B_{2}B_{2}^{\ast }P$ generates by hypotheses an exponentially
stable semigroup on $H$ (see (\ref{2-0}))$.$ The solution to (\ref{26}) is 
\[
y_{N}(t)=e^{\Lambda t}y_{0}+\int_{0}^{t}e^{\Lambda (t-s)}\left\{
B_{2}B_{2}^{\ast }(P-P_{N})y_{N}(s)+B_{1}w(s)\right\} ds 
\]%
which implies that%
\[
\left\Vert y_{N}(t)\right\Vert _{H}\leq e^{-\alpha t}\left\Vert
y_{0}\right\Vert _{H}+C\int_{0}^{t}e^{-\alpha (t-s)}\left( \left\Vert
(P-P_{N})y_{N}(s)\right\Vert _{H}+\left\Vert w(s)\right\Vert _{W}\right) ds, 
\]%
where we took into account the continuity of operators $B_{1}$ and $B_{2}$.
Moreover, we have 
\begin{equation}
\left\Vert y_{N}(t)\right\Vert _{H}^{2}\leq C\left\{ ((f\ast
g)(t))^{2}+\left\Vert y_{0}\right\Vert _{H}^{2}e^{-2\alpha t}\right\} ,
\label{26-00}
\end{equation}%
where $f(t)=e^{-\alpha t}$ and $g(t)=\left\Vert (P-P_{N})y_{N}(t)\right\Vert
_{H}+\left\Vert w(t)\right\Vert _{W}.$ Recall that $C$ takes different
values from a line to another. By recalling the Young inequality for
convolutions 
\[
\left\Vert f\ast g\right\Vert _{L^{r}(0,\infty )}\leq \left\Vert
f\right\Vert _{L^{p}(0,\infty )}\left\Vert g\right\Vert _{L^{q}(0,\infty )},%
\mbox{ with }\frac{1}{p}+\frac{1}{q}=1+\frac{1}{r} 
\]%
and applying it to (\ref{26-00}) for $r=2,$ $p=1$ and $q=2,$ we get 
\begin{eqnarray*}
&&\int_{0}^{\infty }\left\Vert y_{N}(t)\right\Vert _{H}^{2}\leq
C\int_{0}^{\infty }((f\ast g)(t))^{2}dt+C\left\Vert y_{0}\right\Vert
_{H}^{2}\int_{0}^{\infty }e^{-2\alpha s}ds \\
&\leq &C\left( \int_{0}^{\infty }e^{-\alpha t}dt\left( \int_{0}^{\infty
}g^{2}(t)dt\right) ^{1/2}\right) ^{2}+C\left\Vert y_{0}\right\Vert _{H}^{2}
\\
&\leq &C\left( \int_{0}^{\infty }e^{-\alpha t}dt\right) ^{2}\int_{0}^{\infty
}g^{2}(t)dt+C\left\Vert y_{0}\right\Vert _{H}^{2} \\
&\leq &C\int_{0}^{\infty }\left( \left\Vert (P-P_{N})y_{N}(t)\right\Vert
_{H}^{2}+\left\Vert w(t)\right\Vert _{W}^{2}\right) dt+C\left\Vert
y_{0}\right\Vert _{H}^{2} \\
&\leq &C\int_{0}^{\infty }\left\Vert w(t)\right\Vert
_{W}^{2}dt+C\int_{0}^{\infty }\left\Vert (P-P_{N})y_{N}(t)\right\Vert
_{H}^{2}dt+C\left\Vert y_{0}\right\Vert _{H}^{2},
\end{eqnarray*}%
hence 
\begin{equation}
\int_{0}^{\infty }\left\Vert y_{N}(t)\right\Vert _{H}^{2}\leq
C_{0}+C\int_{0}^{\infty }\left\Vert (P-P_{N})y_{N}(t)\right\Vert _{H}^{2}dt,
\label{26-0}
\end{equation}%
due to the fact that $w\in L^{2}(0,\infty ;W).$ Here $C_{0}$ is another
positive constant.

Now, we estimate the second term. We have, by the definition of $P_{N}$ and
the Parseval inequality%
\[
\left\Vert (P-P_{N})\widetilde{y}\right\Vert _{H}^{2}=\left\Vert
\sum_{j=N+1}^{\infty }\gamma _{j}(\widetilde{y},e_{j})_{H}e_{j}\right\Vert
_{H}^{2}=\left\Vert \widetilde{y}\right\Vert _{H}^{2}\sum_{j=N+1}^{\infty
}\gamma _{j}^{2},\mbox{ for }\widetilde{y}\in H. 
\]%
We recall (\ref{52-3}) and take into account that, since $P\in \mathcal{H}$
by Theorem \ref{T1.1}, the series 
\begin{equation}
\left\Vert P\right\Vert _{\mathcal{H}}^{2}=\sum_{j=1}^{\infty }\gamma
_{j}^{2}  \label{26-1}
\end{equation}%
is convergent. Then, its rest, which depends on $N,$ tends to $0,$ that is 
\begin{equation}
\sum_{j=N+1}^{\infty }\gamma _{j}^{2}=:\varepsilon _{N},\mbox{ }\varepsilon
_{N}\rightarrow 0\mbox{ as }N\rightarrow \infty .  \label{27}
\end{equation}%
Thus,%
\begin{equation}
\left\Vert (P-P_{N})\widetilde{y}\right\Vert _{H}^{2}\leq \varepsilon
_{N}\left\Vert \widetilde{y}\right\Vert _{H}^{2},\mbox{ for all }\widetilde{y%
}\in H.  \label{28}
\end{equation}%
We go back to (\ref{26-0}) and deduce that 
\[
\int_{0}^{\infty }\left\Vert y_{N}(t)\right\Vert _{H}^{2}\leq
C_{0}+C\varepsilon _{N}\int_{0}^{\infty }\left\Vert y_{N}(t)\right\Vert
_{H}^{2}dt 
\]%
which yields%
\begin{equation}
(1-C\varepsilon _{N})\int_{0}^{\infty }\left\Vert y_{N}(t)\right\Vert
_{H}^{2}dt\leq C_{0}.  \label{28-0}
\end{equation}%
This allows us to conclude that for $N$ sufficiently large ($\varepsilon
_{N} $ sufficiently small) we get (\ref{25}). By the Datko's result (see 
\cite{Datko-72}) this is equivalent with the fact that the semigroup $%
\Lambda _{N}=-A-B_{2}B_{2}^{\ast }P_{N}$ is exponentially stable.

\noindent To obtain (\ref{24}) we multiply (\ref{20}) by $y_{N}(t)$ and
integrate from $0$ to $t.$ We have 
\begin{eqnarray}
&&\frac{1}{2}\left\Vert y_{N}(t)\right\Vert _{H}^{2}+\omega
\int_{0}^{t}\left\Vert y_{N}(\tau )\right\Vert _{V}^{2}d\tau  \label{29} \\
&\leq &\frac{1}{2}\left\Vert y_{0}\right\Vert
_{H}^{2}+\int_{0}^{t}\left\Vert B_{1}w(s)\right\Vert _{H}\left\Vert
y_{N}(s)\right\Vert _{H}ds+\int_{0}^{t}\left\Vert B_{2}B_{2}^{\ast
}P_{N}y_{N}(s)\right\Vert _{H}\left\Vert y_{N}(s)\right\Vert _{H}ds 
\nonumber \\
&\leq &\frac{1}{2}\left\Vert y_{0}\right\Vert _{H}^{2}+C\int_{0}^{t}\left(
\left\Vert w(s)\right\Vert _{W}\left\Vert y_{N}(s)\right\Vert
_{H}+\left\Vert y_{N}(s)\right\Vert _{H}^{2}\right) ds  \nonumber \\
&\leq &C\int_{0}^{t}\left\Vert y_{N}(s)\right\Vert _{H}^{2}ds+C\left(
\left\Vert y_{0}\right\Vert _{H}^{2}+\int_{0}^{t}\left\Vert w(s)\right\Vert
_{W}^{2}ds\right) .  \nonumber
\end{eqnarray}%
Here, we used the continuity of all linear operators $B_{1},$ $B_{2},$ $%
B_{2}^{\ast }$ and (\ref{30}). Hence, by (\ref{29}) we obtain 
\begin{equation}
\left\Vert y_{N}(t)\right\Vert _{H}^{2}\leq C\left( \left\Vert
y_{0}\right\Vert _{H}^{2}+\int_{0}^{t}\left\Vert w(s)\right\Vert
_{W}^{2}ds\right) e^{Ct},\mbox{ for all }t>0,  \label{31}
\end{equation}%
that is (\ref{24}). Also, (\ref{29}) implies $\left\Vert y\right\Vert
_{L^{2}(0,T;V)}\leq C$.

Further, let us subtract equations satisfied by $y$ and $y_{N}$ written in
form (\ref{26}), 
\begin{equation}
(y(t)-y_{N}(t))^{\prime }=(-A-B_{2}B_{2}^{\ast
}P)(y(t)-y_{N}(t))-B_{2}B_{2}^{\ast }(P-P_{N})y_{N}(t),  \label{32}
\end{equation}%
and write its solution 
\begin{equation}
y(t)-y_{N}(t)=\int_{0}^{t}e^{\Lambda (t-s)}(-B_{2}B_{2}^{\ast
}(P-P_{N})y_{N}(s))ds,  \label{80}
\end{equation}%
implying 
\[
\left\Vert y(t)-y_{N}(t)\right\Vert _{H}\leq \int_{0}^{t}e^{-\alpha
(t-s)}\left\Vert (P-P_{N})y_{N}(s)\right\Vert _{H}ds. 
\]%
Next, by a similar calculus as before, using the Young inequality for
convolutions for $f(t)=e^{-\alpha t}$ and $g(t)=\left\Vert
(P-P_{N})y_{N}(s)\right\Vert _{H},$ and (\ref{28}) we have 
\begin{eqnarray}
&&\int_{0}^{\infty }\left\Vert y(t)-y_{N}(t)\right\Vert _{H}^{2}dt\leq
C\left( \int_{0}^{t}e^{-\alpha (t-s)}\left\Vert (P-P_{N})y_{N}(s)\right\Vert
_{H}ds\right) ^{2}  \label{81} \\
&\leq &C\int_{0}^{\infty }\left\Vert (P-P_{N})y_{N}(t)\right\Vert
_{H}^{2}dt\leq C\varepsilon _{N}\int_{0}^{\infty }\left\Vert
y_{N}(t)\right\Vert _{H}^{2}dt.  \nonumber
\end{eqnarray}%
Recalling (\ref{25}) we obtain that for $N$ large 
\begin{equation}
\int_{0}^{\infty }\left\Vert y(t)-y_{N}(t)\right\Vert _{H}^{2}dt\leq
C\varepsilon _{N}.  \label{82}
\end{equation}%
By multiplying (\ref{32}) by $(y-y_{N})(t)$ and integrating from $0$ to $t$
we obtain after a few calculations 
\begin{eqnarray*}
&&\frac{1}{2}\left\Vert y(t)-y_{N}(t)\right\Vert _{H}^{2}+\omega
\int_{0}^{t}\left\Vert y(s)-y_{N}(s)\right\Vert _{H}^{2}ds \\
&\leq &\int_{0}^{t}(-B_{2}B_{2}^{\ast }P(y(s)-y_{N}(s)),y(s)-y_{N}(s))_{H}ds
\\
&&+\int_{0}^{t}\left( -B_{2}B_{2}^{\ast
}(P-P_{N})y_{N}(s),y(s)-y_{N}(s)\right) _{H}ds \\
&\leq &C\int_{0}^{t}\left\Vert y(s)-y_{N}(s)\right\Vert
_{H}^{2}ds+C\int_{0}^{t}\left\Vert (P-P_{N})y_{N}(s)\right\Vert _{H}^{2}ds \\
&\leq &C\int_{0}^{t}\left\Vert y(s)-y_{N}(s)\right\Vert
_{H}^{2}ds+C\varepsilon _{N}\int_{0}^{t}\left\Vert y_{N}(s)\right\Vert
_{H}^{2}ds,
\end{eqnarray*}%
whence, by (\ref{25})%
\begin{equation}
\left\Vert y(t)-y_{N}(t)\right\Vert _{H}^{2}\leq C\varepsilon _{N}e^{Ct}%
\mbox{ for all }t>0.  \label{83}
\end{equation}%
By (\ref{28}) we get 
\begin{eqnarray}
&&\left\Vert (u-u_{N})(t)\right\Vert _{U}^{2}=\left\Vert B_{2}^{\ast }\left(
Py(t)-P_{N}y_{N}(t)\right) \right\Vert _{U}^{2}  \label{35} \\
&\leq &C\left\Vert (P-P_{N})y(t)\right\Vert _{H}^{2}+C\left\Vert
P_{N}(y(t)-y_{N}(t))\right\Vert _{H}^{2}  \nonumber \\
&\leq &C\varepsilon _{N}\left\Vert y(t)\right\Vert _{H}^{2}+C\left\Vert
y(t)-y_{N}(t)\right\Vert _{H}^{2}\left\Vert P\right\Vert _{\mathcal{H}}^{2},%
\mbox{ for all }t>0,\mbox{ }y\in C([0,\infty );H).  \nonumber
\end{eqnarray}%
Moreover, by (\ref{82}) we still have 
\begin{equation}
\int_{0}^{\infty }\left\Vert (u-u_{N})(t)\right\Vert _{U}^{2}dt\leq
C\varepsilon _{N}\int_{0}^{\infty }\left\Vert y(t)\right\Vert
_{H}^{2}dt+C\int_{0}^{\infty }\left\Vert y(t)-y_{N}(t)\right\Vert
_{H}^{2}dt\leq C_{1}\varepsilon _{N}.  \label{35-0}
\end{equation}%
Finally, using (\ref{83}) and (\ref{35}) we have 
\begin{eqnarray}
&&\left\Vert z_{N}(t)\right\Vert _{Z}=\left\Vert
C_{1}y_{N}(t)+D_{1}u_{N}(t)\right\Vert _{Z}  \nonumber \\
&=&\left\Vert
C_{1}(y_{N}(t)-y(t))+D_{1}(u_{N}(t)-u(t))+C_{1}y(t)+D_{1}u(t)\right\Vert _{Z}
\nonumber \\
&\leq &\left\Vert C_{1}(y_{N}(t)-y(t))\right\Vert _{Z}+\left\Vert
D_{1}(u_{N}(t)-u(t))\right\Vert _{Z}+\left\Vert z(t)\right\Vert _{Z}
\label{36} \\
&\leq &C\left( \left\Vert y_{N}(t)-y(t)\right\Vert _{H}+\left\Vert
u_{N}(t)-u(t)\right\Vert _{U}\right) +\left\Vert z(t)\right\Vert _{Z}. 
\nonumber
\end{eqnarray}%
Therefore, (\ref{82}) and (\ref{35-0}) yield 
\begin{equation}
\left\Vert z_{N}\right\Vert _{L^{2}(\mathbb{R}_{+}:Z)}^{2}=\int_{0}^{\infty
}\left\Vert z_{N}(t)\right\Vert _{Z}^{2}dt\leq C\varepsilon
_{N}+\int_{0}^{\infty }\left\Vert z(t)\right\Vert _{Z}^{2}dt\leq
C\varepsilon _{N}+\gamma ^{2}\left\Vert w\right\Vert _{L^{2}(\mathbb{R}%
_{+}:W)}^{2}  \label{36-0}
\end{equation}%
proving that $z_{N}$ satisfies relation (\ref{22}) with $\delta
_{N}=\varepsilon _{N}^{1/2}\rightarrow 0,$ as $N\rightarrow \infty .$ In
conclusion we have proved that $y_{N}$ is the unique solution to (\ref{20})
satisfying (\ref{24})-(\ref{25-00}). Eventually, since $z_{N}$ satisfies (%
\ref{36-0}) it follows that $u_{N}=-B_{2}^{\ast }P_{N}y_{N}$ is the
approximating solution to the $H^{\infty }$-control problem associated to
system (\ref{20})-(\ref{21}), satisfying (\ref{35})-(\ref{35-0}).

Now, we can pass to the limit as $N\rightarrow \infty .$ By (\ref{24}) and (%
\ref{83}) we deduce that 
\[
y_{N}(t)\rightarrow y(t)\mbox{ strongly in }H,\mbox{ uniformly on }[0,T],%
\mbox{ for all }T>0, 
\]%
while by (\ref{82})%
\[
y_{N}\rightarrow y\mbox{ strongly in }L^{2}(0,\infty ;H). 
\]%
By (\ref{35}) and (\ref{35-0}) it follows that $u_{N}(t)\rightarrow u(t)$
strongly in $U$ a.e. $t>0$ and 
\[
u_{N}\rightarrow u\mbox{ strongly in }L^{2}(0,\infty ;U).\mbox{ } 
\]%
By (\ref{36}) and (\ref{36-0}) it follows that $z_{N}(t)\rightarrow z(t)$
strongly in $Z$ a.e. $t>0$ and 
\[
z_{N}\rightarrow z\mbox{ strongly in }L^{2}(0,\infty ;Z).\mbox{ } 
\]%
Then, by (\ref{25-1}) we get 
\[
P_{N}y\rightarrow Py\mbox{ strongly in }H,\mbox{ for all }y\in H 
\]%
and by (\ref{25-0}), 
\[
P_{N}\rightarrow P\mbox{ weakly in }\mathcal{H}\mbox{ as }N\rightarrow
\infty . 
\]%
All these show that the $u_{N}$ suitably approximates the feedback control
solution to the original $H^{\infty }$-control problem.\hfill $\square $

\section{The noncoercive case}

\setcounter{equation}{0}

In the noncoercive case, Theorem \ref{T1.2} states the existence of a unique
solution to the Riccati equation and of the $H^{\infty }$-control problem.
The aspect arose by the noncoercivity of $A$ in proving the existence of the
solution of the Riccati equation is overcome by considering a coercive
approximation of $A$ by the operator $A_{\omega }=A+\omega J,$ with $J$ the
canonical isomorphism from $V$ to $V^{\prime }$ (see Theorem 3.1 in \cite%
{GM-2026}). The solution of Riccati equation $P_{\omega }$ exists and it is
unique and then it is shown that its weak limit in $\mathcal{V},$ as $\omega
\rightarrow 0,$ solves equation (\ref{5}) with $A\geq 0.$ Thus, this case
reduces to approaching the $H^{\infty }$-control problem for the system (\ref%
{1})-(\ref{2}) with $A$ replaced by $A_{\omega }.$ To calculate the
truncated solution, Theorem \ref{T2.2} will be applied to equations (\ref{20}%
)-(\ref{21}) with $A$ replaced by $A_{\omega }$ and this yields a solution $%
P_{\omega }^{N}$ given by%
\[
P_{\omega }^{N}y=\sum_{j=1}^{N}\gamma _{j}^{\omega }y_{j}e_{j} 
\]%
where $\gamma _{j}^{\omega }$ has the expression 
\begin{equation}
\gamma _{j}^{\omega }=\frac{-a_{jj}^{\omega }+\sqrt{(a_{jj}^{\omega
})^{2}+g_{jj}\left\Vert C_{1}e_{j}\right\Vert _{H}^{2}}}{g_{jj}},\mbox{ }
\label{55-7}
\end{equation}%
with 
\[
a_{jj}^{\omega }:=\left\langle A_{\omega }e_{j},e_{j}\right\rangle
_{V^{\prime },V},\mbox{ }g_{jj}:=(\Gamma e_{j},e_{j})_{H},\mbox{ }j\geq 1. 
\]%
Here, $a_{jj}^{\omega }\geq \omega $ and $g_{jj}\geq g_{0}>0$ by hypothesis (%
\ref{52}), implying that 
\begin{equation}
0<\gamma _{j}^{\omega }\leq \frac{\left\Vert C_{1}e_{j}\right\Vert _{H}}{%
\sqrt{g_{jj}}}\leq \frac{\left\Vert C_{1}e_{j}\right\Vert _{H}}{\sqrt{g_{0}}}%
\leq \frac{1}{\sqrt{g_{0}}}\left\Vert C_{1}\right\Vert _{\mathcal{H}}.
\label{55-6}
\end{equation}%
Then, Theorem \ref{T2.2} can be applied for the system with the operator $%
A_{\omega }$ to get an approximating solution $P_{\omega }^{N}.$ Moreover, $%
P_{\omega }^{N}\rightarrow P_{\omega }$ as $N\rightarrow \infty $ as in
Theorem \ref{T2.2}.

\section{An alternative approximating system}

\setcounter{equation}{0}

In the case when the operators $C_{1}$ and $\Gamma $ are not Hilbert-Schmidt
and the existence of the Riccati equation cannot not be proved by Theorems %
\ref{T1.1} and \ref{T1.2}, we provide in this section another way of
approximating the system to be controlled. Thus, considering the orthonormal
bases $\{e_{j}\}_{j\geq 1},$ $\{e_{j}^{u}\}_{j\geq 1},$ $\{e_{j}^{w}\}_{j%
\geq 1},$ $\{e_{j}^{z}\}_{j\geq 1},$ in the spaces $H,$ $U,$ $W,$ $Z,$
respectively we define the truncated series 
\begin{equation}
\sum_{j=1}^{N}y_{j}e_{j},\mbox{ }\sum_{j=1}^{N}u_{j}e_{j}^{u},\mbox{ }%
\sum_{j=1}^{N}w_{j}e_{j}^{w},\mbox{ }\sum_{j=1}^{N}z_{j}e_{j}^{z},
\label{60}
\end{equation}%
introduce them in equations (\ref{1}) and (\ref{2}) and obtain%
\[
\sum_{j=1}^{N}y_{j}^{\prime
}e_{j}+\sum_{j=1}^{N}y_{j}Ae_{j}=\sum_{j=1}^{N}u_{j}B_{2}e_{j}^{u}+%
\sum_{j=1}^{N}w_{j}B_{1}e_{j}^{w}, 
\]%
\[
\sum_{j=1}^{N}z_{j}e_{j}^{z}=\sum_{j=1}^{N}y_{j}C_{1}e_{j}+%
\sum_{j=1}^{N}u_{j}D_{1}e_{j}^{u}. 
\]%
By multiplying by $e_{i}$ the first equation and by $e_{i}^{z}$ the second
one, we get 
\begin{equation}
y_{i}^{\prime
}+\sum_{j=1}^{N}a_{ij}y_{j}=\sum_{j=1}^{N}u_{j}b_{ij}^{2}+%
\sum_{j=1}^{N}w_{j}b_{ij}^{1},\mbox{ }i=1,...,N,  \label{61}
\end{equation}%
\begin{equation}
z_{i}=\sum_{j=1}^{N}y_{j}c_{ij}^{1}+\sum_{j=1}^{N}u_{j}d_{ij}^{1},\mbox{ }%
i=1,...,N,  \label{61-0}
\end{equation}%
where 
\begin{equation}
a_{ij}:=\left\langle A^{\ast }e_{i},e_{j}\right\rangle _{V^{\prime },V},%
\mbox{ }b_{ij}^{2}:=\left( B_{2}^{\ast }e_{i},e_{j}^{u}\right) _{U},\mbox{ }%
b_{ij}^{1}:=\left( B_{1}^{\ast }e_{i},e_{j}^{w}\right) _{W},\mbox{ }
\label{62}
\end{equation}%
\[
c_{ij}^{1}:=\left( C_{1}^{\ast }e_{i}^{z},e_{j}\right) _{H},\mbox{ }%
d_{ij}^{1}:=\left( D_{1}^{\ast }e_{i}^{z},e_{j}^{u}\right) _{U}. 
\]%
We introduce the column vectors 
\[
y^{N}=(y_{i})_{i=1,...N},\mbox{ }u^{N}=(u_{i})_{i=1,...N},\mbox{ }%
w^{N}=(w_{i})_{i=1,...N},\mbox{ }z^{N}=(z_{i})_{i=1,...N}, 
\]%
the $N\times N$ matrices%
\[
A_{N}=(a_{ij})_{i,j=1,...N},\mbox{ }B_{2N}=(b_{ij}^{2})_{i,j=1,...N},\mbox{ }%
B_{1N}=(b_{ij}^{1})_{i,j=1,...N}, 
\]%
\[
C_{1N}=(c_{ij}^{1})_{i,j=1,...N},\mbox{ }D_{1N}=(d_{ij}^{1})_{i,j=1,...N} 
\]%
and write (\ref{61})-(\ref{61-0}) as 
\begin{equation}
(y^{N})^{\prime }+A_{N}y^{N}=B_{2N}u^{N}+B_{1N}w^{N},  \label{63}
\end{equation}%
with 
\begin{equation}
z^{N}=C_{1N}y^{N}+D_{1N}u^{N}.  \label{63-0}
\end{equation}%
Here, $A_{N}$ preserves the properties of $A,$ that is it is positive if $A$
is positive. The system was reduced to the finite dimensional case where the
spaces $H,$ $V$ and $V^{\prime }$ are actually $\mathbb{R}^{N}.$ In this
case, $B_{2N}$ and $B_{1N}$ are Hilbert-Schmidt, such that Theorem \ref{T1.1}
can be applied to show that the associated Riccati equation 
\begin{equation}
A_{N}^{\ast }P_{N}+P_{N}A_{N}+P_{N}(B_{2N}B_{2N}^{\ast }-\gamma
^{-2}B_{1N}B_{1N}^{\ast })P_{N}=C_{1N}^{\ast }C_{1N}  \label{64}
\end{equation}%
has a unique solution $P_{N},$ and the feedback controller is $%
u_{N}=-B_{2N}^{\ast }P_{N}y_{N}.$ The convergence of $y^{N}$ to $y$ the
solution to (\ref{1}) and of $u_{N}$ to $u$ can be proved in a similar way
as in Theorem \ref{T2.2}.

\section{An example}

\setcounter{equation}{0}

We present in this section the construction of the approximating solutions
to the Riccati equation and to the $H^{\infty }$-control problem for the
system%
\begin{equation}
y_{t}-\Delta y=B_{1}w+B_{2}u,\mbox{ in }(0,\infty )\times \Omega ,
\label{21-0}
\end{equation}%
\begin{equation}
y=0,\mbox{ on }(0,\infty )\times \partial \Omega ,  \label{21-1}
\end{equation}%
\begin{equation}
y(0)=0,\mbox{ in }\Omega ,  \label{21-2}
\end{equation}%
\begin{equation}
z=C_{1}y,\mbox{ in }(0,\infty )\times \Omega ,  \label{21-3}
\end{equation}%
in the domain $\Omega =(0,1)$. Let us consider the variational triplet $%
H_{0}^{1}(\Omega )\subset L^{2}(\Omega )\subset H^{-1}(\Omega ),$ define the
operator 
\[
-\Delta :H_{0}^{1}(\Omega )\rightarrow H^{-1}(\Omega ),\mbox{ }\left\langle
-\Delta y,\psi \right\rangle _{H^{-1}(\Omega ),H_{0}^{1}(\Omega
)}=\int_{\Omega }\nabla y\cdot \nabla \psi dx, 
\]%
and consider, as example, the operators 
\begin{eqnarray}
B_{1} &:&L^{2}(\Omega )\rightarrow H_{0}^{1}(\Omega ),\mbox{ }%
B_{1}w=(-\Delta )^{-1}w,\mbox{ \ }  \label{21-7} \\
B_{2} &:&L^{2}(\Omega )\rightarrow L^{2}(\Omega ),\mbox{ }%
(B_{2}u)(x)=\int_{\Omega }k(x,\xi )u(\xi )d\xi ,  \nonumber \\
C_{1} &:&L^{2}(\Omega )\rightarrow H_{0}^{1}(\Omega ),\mbox{ }%
C_{1}y=(-\Delta )^{-1}y,  \nonumber
\end{eqnarray}%
where 
\begin{equation}
k\in L^{2}(\Omega \times \Omega ),\mbox{ }k(x,\xi )=\exp (-\left\vert x-\xi
\right\vert ^{2}).  \label{21-7-0}
\end{equation}%
We set the spaces 
\begin{equation}
H=W=Z=L^{2}(\Omega ),\mbox{ }V=H_{0}^{1}(\Omega ),\mbox{ }U=L^{2}(\Omega )
\label{21-6}
\end{equation}%
and define the adjoint operators 
\begin{eqnarray}
B_{1}^{\ast } &:&H\rightarrow V,\mbox{ }B_{1}^{\ast }v=(-\Delta )^{-1}v,%
\mbox{ for all }v\in H,  \label{22-3} \\
B_{2}^{\ast } &:&H\rightarrow H,\mbox{ }(B_{2}^{\ast }v)(x)=\int_{\Omega
}k(\xi ,x)v(\xi )d\xi ,\mbox{ for all }v\in H,  \nonumber \\
C_{1}^{\ast } &:&H\rightarrow V,\mbox{ }C_{1}^{\ast }v=(-\Delta )^{-1}v,%
\mbox{ for all }v\in H.  \nonumber
\end{eqnarray}%
We introduce the self-adjoint operator $A:D(A)\subset L^{2}(\Omega
)\rightarrow L^{2}(\Omega ),\mbox{ }Ay=-\Delta y,$ with $D(A)=\{y\in
H_{0}^{1}(\Omega );$ $Ay\in L^{2}(\Omega )\}.$ Then, equation (\ref{21-0})
can be equivalently written%
\begin{equation}
y^{\prime }(t)+Ay(t)=B_{1}w(t)+B_{2}u(t),\mbox{ }t\geq 0.  \label{24-1}
\end{equation}%
In order to apply Theorem \ref{T2.2} to get the existence of the Riccati
equation (\ref{5}) we have to check assumptions of Theorem \ref{T1.1}.
Indeed, the operator $A$ is continuous and coercive, the operators $C_{1}$
and $\Gamma $ are Hilbert-Schmidt by a direct calculation of their norms,
and $\Gamma \geq 0$ under some conditions as proved in \cite{GM-2026}, Lemma
5.1. Concerning the latter, we shall give a few explanations. First, we take
as basis $\{e_{j}\}_{j\geq 1}$ for the Fourier expansion of the solution the
eigenvectors of the Laplacian, 
\[
-\Delta e_{j}=\lambda _{j}e_{j},\mbox{ }e_{j}(0)=e_{j}(1)=0,\mbox{ }j\geq 1 
\]%
and recall the Weyl's law characterizing the asymptotic distribution of the
eigenvalues for the Laplacian on a bounded domain, 
\begin{equation}
c_{1}^{\Omega }j^{2/d}<\lambda _{j},\mbox{ for all }j\geq 1,  \label{59}
\end{equation}%
(see e.g., \cite{Li-Yau-1983}), where $c_{1}^{\Omega }$ is a positive
constant depending on $\Omega ,$ and $\,d$ is the space dimension (here $%
d=1).$ Then, the property $\Gamma \geq 0$ will follow, as proved in \cite%
{GM-2026}, Lemma 5.1, under the following assumptions:%
\begin{equation}
\sum\limits_{j=1}^{\infty }j^{-4/d}\leq \gamma ^{2}(c_{1}^{\Omega
})^{2}\lambda _{1}^{B},  \label{64-1}
\end{equation}%
\begin{equation}
B_{2}B_{2}^{\ast }>0,\mbox{ }  \label{64-0}
\end{equation}%
where $\lambda _{1}^{B}$ is the first eigenvalue of the operator $%
B_{2}B_{2}^{\ast }$ and $\gamma $ is the robust performance index in (\ref{3}%
). We recall that $B_{2}B_{2}^{\ast }>0$ means $(B_{2}B_{2}^{\ast
}y,y)_{H}>0 $ for all $y\in H,$ $y\neq 0.$

Now, we check if these conditions are satisfied in our example.

The first condition (\ref{64-1}) can be obtained by choosing a sufficiently
large $\gamma ,$ as we shall see a little later. The second one, (\ref{64-0}%
), means that $\left( B_{2}B_{2}^{\ast }v,v\right) _{H}=\left\Vert
B_{2}^{\ast }v\right\Vert _{H}^{2}>0,$ for all $v\in H,$ that is this norm
cannot vanish. This amounts to say that 
\begin{equation}
(B_{2}^{\ast }v)(x)=0\mbox{ a.e. }x\in \Omega ,\mbox{ }\Longrightarrow v(x)=0%
\mbox{ a.e. }x\in \Omega \subset \mathbb{R}^{d}.  \label{70}
\end{equation}%
This is proved in the lemma below, for a kernel $k$ of a more general form,
satisfying certain properties, which in particular are valid also for the
kernel (\ref{21-7-0}) in our example.

\begin{lemma}
\label{L5.1}Assume that $k(x,x^{\prime })=p(x-x^{\prime }),$ for $x,$ $%
x^{\prime }\in \mathbb{R}^{d},$ where $p\in L^{1}(\mathbb{R}),$ is analytic
and 
\begin{equation}
\widehat{p}(\xi )\neq 0,\mbox{ for all }\xi \in \mathbb{R}^{d},  \label{71}
\end{equation}%
where $\widehat{p}(\xi )=\int_{\mathbb{R}^{d}}e^{-i\xi x}p(x)dx.$ Then, (\ref%
{70}) holds.
\end{lemma}

\noindent \textbf{Proof. }Let $v\in L^{2}(\Omega ),$ let $\varphi
=B_{2}^{\ast }v,$ and so 
\begin{equation}
\varphi (x)=\int_{\Omega }p(x^{\prime }-x)v(x^{\prime })dx^{\prime }=\int_{%
\mathbb{R}^{d}}p(x^{\prime }-x)\widetilde{v}(x^{\prime })dx^{\prime },\mbox{ 
}x\in \mathbb{R}^{d},  \label{72}
\end{equation}%
where $\widetilde{v}(x)=v\left( x\right) $ for $x\in \Omega $ and $%
\widetilde{v}(x)=0$ for $x\in \mathbb{R}^{d}\backslash \Omega .$

Then, $\varphi =p\ast \widetilde{v}$ and we have 
\begin{equation}
\widehat{\varphi }(\xi )=\widehat{p}(\xi )\widehat{\widetilde{v}}(\xi )\mbox{
for }\xi \in \mathbb{R}^{d}.  \label{73}
\end{equation}%
If $(B_{2}^{\ast }v)(x)=0$ we have $\varphi (x)=0$ for all $x\in \mathbb{R}%
^{d},$ but since $\varphi $ is analytic (because $p$ is) it follows that $%
\widehat{\varphi }(\xi )=0$ for all $\xi \in \mathbb{R}^{d}.$ Then, by (\ref%
{71}) and (\ref{73}) it follows that $\widehat{\widetilde{v}}(\xi )=0$ for
all $\xi \in \mathbb{R}^{d}.$ Since by Parseval formula%
\[
(2\pi )^{d}\int_{\mathbb{R}^{N}}(\widetilde{v}(x))^{2}dx=\int_{\mathbb{R}%
^{N}}(\widehat{\widetilde{v}}(\xi ))^{2}d\xi 
\]%
we get $\widetilde{v}\equiv 0$ and therefore, $v=0$ on $\Omega ,$ as
claimed.\hfill $\square $

\smallskip

Examples are $p(x)=1/\left( \left\Vert x\right\Vert _{\mathbb{R}%
^{d}}^{2}+1\right) ^{m},$ $\alpha >0,$ $m>d/2,$ or $p(x)=\exp (-\alpha
\left\Vert x\right\Vert _{\mathbb{R}^{d}}^{2})$, as it is $k$ in our example.

\smallskip

We return now to our example. Since $C_{1}$ and $\Gamma $ are
Hilbert-Schmidt and $\Gamma \geq 0,$ by the Lemma before, under conditions (%
\ref{64-1}) and (\ref{64-0}), the $H^{\infty }$-control problem for the
system (\ref{21-0})-(\ref{21-3}) has a solution $u=-B_{2}^{\ast }Py$ where $%
P $ is the unique solution to the Riccati equation (\ref{5}). These
solutions are approximated by (\ref{13-0})-(\ref{13}) 
\begin{eqnarray}
P_{N}y_{N} &=&\sum_{j=1}^{N}\gamma _{j}y_{j}e_{j},\mbox{ }u_{N}=F_{N}y_{N},%
\mbox{ }  \label{74} \\
F_{N}y_{N} &=&-B_{2}^{\ast }P_{N}y_{N}=-\sum_{j=1}^{N}\gamma
_{j}y_{j}B_{2}^{\ast }e_{j}=-\sum_{j=1}^{N}\gamma _{j}y_{j}\int_{0}^{1}k(\xi
,x)e_{j}(\xi )d\xi ,  \nonumber
\end{eqnarray}%
where $\{e_{j}\}_{j\geq 1}$ represent the eigenvectors of the operator $%
-\Delta $ previously defined, that is $-\Delta e_{j}=\lambda _{j}e_{j}$ and $%
\gamma _{j}$ are given by (\ref{52-2}).

In this case, $P_{N}$ and the control $F_{N}$ follow to be the diagonal
matrices 
\[
P_{N}=\left( 
\begin{array}{ccc}
\gamma _{1}e_{1} &  &  \\ 
& ... &  \\ 
&  & \gamma _{N}e_{N}%
\end{array}%
\right) ,\mbox{ } 
\]%
\[
F_{N}=\left( 
\begin{array}{ccc}
-\gamma _{1}\int_{\Omega }k(\xi ,x)e_{1}(\xi )d\xi &  &  \\ 
& ... &  \\ 
&  & -\gamma _{N}\int_{\Omega }k(\xi ,x)e_{N}(\xi )d\xi%
\end{array}%
\right) . 
\]%
Finally, we make some comments about $\gamma _{j}.$ In this case the
eigenvectors and eigenvalues are 
\[
e_{j}(x)=\sqrt{2}\sin \pi jx,\mbox{ }\lambda _{j}=j^{2}\pi ^{2},\mbox{ }%
j\geq 1, 
\]%
\[
a_{kj}=\left\langle Ae_{k},e_{j}\right\rangle _{V^{\prime },V}=\left\langle
-\Delta e_{k},e_{j}\right\rangle _{V^{\prime }V}=\lambda _{k}\delta _{kj}, 
\]%
and since $k(x,\xi )=\exp (-\left\vert x-\xi \right\vert ^{2}),$ we have 
\begin{eqnarray*}
g_{kj} &=&(\Gamma e_{k},e_{j})_{H}=((B_{2}B_{2}^{\ast }-\gamma
^{-2}B_{1}B_{1}^{\ast }e_{k},e_{j})_{H} \\
&=&(B_{2}^{\ast }e_{k},B_{2}^{\ast }e_{j})_{H}-\gamma ^{-2}(B_{1}^{\ast
}e_{k},B_{1}^{\ast }e_{j})_{H}.
\end{eqnarray*}%
Denoting 
\[
I_{j}(x)=B_{2}^{\ast }e_{j}(x)=\sqrt{2}\int_{\Omega }\exp (-\left\vert \xi
-x\right\vert ^{2})\sin \pi j\xi d\xi 
\]%
and 
\[
I_{kj}=\left( I_{k},I_{j}\right) _{H} 
\]%
we obtain 
\begin{eqnarray*}
g_{kj} &=&\int_{0}^{1}I_{k}(x)I_{j}(x)dx-\gamma ^{-2}((-\Delta
)^{-1}e_{k},(-\Delta )^{-1}e_{j})_{H} \\
&=&I_{k}{}_{j}-\frac{\gamma ^{-2}}{\lambda _{k}\lambda _{j}}\delta _{kj},%
\mbox{ if }k+j\mbox{ is even and }g_{kj}=0\mbox{ if }k+j\mbox{ is odd.}
\end{eqnarray*}%
By Lemma \ref{L5.1} it follows that $\left\Vert B_{2}^{\ast
}e_{j}\right\Vert _{H}>0$ for all $j\geq 1$ (otherwise $e_{j}$ would be $0)$
and so, this ensures that $I_{kj}>0.$ The, $I_{j}$ behaves 
\[
I_{j}(x)\approx \frac{C}{j}\mbox{ as }j\rightarrow \infty 
\]%
and so 
\[
I_{kj}\approx \frac{C}{kj} 
\]%
and 
\[
g_{kj}\approx \frac{C}{\pi ^{2}kj}-\frac{\gamma ^{-2}}{\pi ^{4}k^{2}j^{2}}%
\mbox{ if }k+j\mbox{ is even and }g_{kj}=0\mbox{ if }k+j\mbox{ is odd.} 
\]%
Thus, the condition for having $g_{jj}\geq 0$ is obtained for $\gamma $
large enough, that is $\gamma ^{2}\succeq \frac{1}{I_{11}\pi ^{4}}.$
Recalling the other condition (\ref{64-1}), where $\sum\limits_{j=1}^{\infty
}j^{-4}=\pi ^{4}/90,$ we deduce that it is sufficient to take 
\[
\gamma ^{2}\geq \max \left\{ \frac{1}{I_{11}\pi ^{4}},\frac{\pi ^{4}}{90}%
\frac{1}{(c_{1}^{\Omega })^{2}\lambda _{1}^{B}}\right\} 
\]%
in order to satisfy the requirement $\Gamma \geq 0.$ In conclusion, $\gamma
_{j}$ can be calculated by (\ref{52-2}) and used to construct the
approximation $P_{N}y_{N}$ by (\ref{74}).

\textbf{Acknowledgment}

The author acknowledges the support of the project
PN-IV-P6-6.1-CoEx-2024-0029, Contract 10CoEx/2026 by the Ministry of
Education and Research, National Research Authority, UEFISCDI.

\end{document}